\documentclass{svmult}

\usepackage{makeidx}         
\usepackage{graphicx}        
\usepackage{multicol}        
\usepackage[bottom]{footmisc}

\makeindex             

\begin{document}

\newcommand{\rdike}{\overleftarrow{\partial}} 
\newcommand{\ldike}{\overrightarrow{\partial}} 

\def\ba{\mbox{\boldmath $A$}}
\def\bb{\mbox{\boldmath $B$}}
\def\bPhiike{\mbox{\boldmath $\Phi$}}
\def\bphiike{\mbox{\boldmath $\phi$}}
\def\bDeltaike{\mbox{\boldmath $\Delta$}}
\def\bdeltaike{\mbox{\boldmath $\delta$}}
\newcommand{\dbvike}{{\Delta_{BV}}}
\newcommand{\sbvike}[2]{{({{#1},{#2}})}}
\def\brsike{\delta}
\newcommand{\OmBVike}{{\Omega_{BV}}}
\newcommand{\bOmBVike}{\boldmath{\Omega_{BV}}}

\newcommand{\calEike}{{\cal E}}
\newcommand{\calMike}{{\cal M}}
\newcommand{\calXike}{{\cal X}}
\newcommand{\calVike}{{\cal V}}
\newcommand{\bracketike}[2]{\langle #1\,,#2\rangle}



\title*{Deformation of Batalin-Vilkovisky Structures}
\author{Noriaki IKEDA}
\institute{Department of Mathematical Sciences, 
Ritsumeikan University \\
Kusatsu, Shiga 525-8577, Japan 
\texttt{nikeda@se.ritsumei.ac.jp}
}
%
%
\maketitle

A Batalin-Vilkovisky formalism is most general framework to construct 
consistent quantum field theories.
Its mathematical structure is called {\it a Batalin-Vilkovisky structure}.
First we explain rather mathematical setting of a Batalin-Vilkovisky formalism.
Next, we consider deformation theory of a Batalin-Vilkovisky structure.
Especially, we consider deformation of topological sigma models 
in any dimension,
which is closely related to deformation theories in mathematics, 
including deformation from commutative geometry to noncommutative geometry.
We obtain a series of new nontrivial topological sigma models 
and we find these models have the Batalin-Vilkovisky structures based on
a series of new algebroids.


\section{Introduction}
\label{sec:1}
Topological field theory is a powerful method to analyze geometry 
by means of a quantum field theory.\cite{BBRTike}
A Batalin-Vilkovisky formalism \cite{BVike}
is a most general systematic method to 
treat a consistent quantum field theory.
So it is natural to treat topological field theory by means of 
a Batalin-Vilkovisky formalism.

Deformation theory is one of main topics in mathematics.
On the other hand, 
deformation theory of a quantum field theory is proposed by
\cite{BBHike}\cite{BHike} in the physical context. 
Purposes of this deformation theory are
to find a new gauge theory, or to prove 
a no-go theorem to construct a new gauge theory.
Now we apply deformation theory to topological field theories,
especially to topological sigma models. 
%
Our purpose is to construct many geometries by topological field theories,
classify geometries as topological field theories, 
and unify many deformation theories as a topological field theories.
Moreover we can analyze many deformation theories 
in mathematics as quantum field theories.

In this article, we consider a {\it topological sigma model}.
Let $X$ and $M$ be two manifolds.
We denote $\phi$ a (smooth) map from $X$ to $M$.
%
A {\it sigma model} is a quantum field theory constructed from a map
$\phi$ (and other auxiliary fields).
We analyze structures on $M$ by the structures induced from $X$.
A {\it topological sigma model} is a sigma model
independent of metrics of $X$ and $M$ \cite{W1ike}.

The AKSZ formulation \cite{AKSZike} of the Batalin-Vilkovisky formalism 
is a general framework to construct a topological sigma model
by the Batalin-Vilkovisky method.
This formulation is appropriate to analyze geometry 
by means of a topological sigma model.
In this article, first we generalize the AKSZ formulation to 
a general $n$ dimensional base manifold $X$ and 
a target manifold with general gradings $p$.

Next we discuss deformation of a Batalin-Vilkovisky structure 
with general gradings.
We can construct various geometrical structures on $M$ if 
we consider general $n$ and general grading $p$.
This formulation provides a clear method to 
unify and to classify many geometries as
Batalin-Vilkovisky structures, and 
to analyze them as quantum field theories.
This construction includes many new topological sigma models
as special cases, for examples,
the topological sigma model
with a K\"ahler structure ($A$ model), 
with a complex structure ($B$ model) \cite{W2ike}\cite{AKSZike},
with a symplectic structure \cite{W1ike},
with a Poisson structure (the Poisson sigma model) \cite{IIike}\cite{SSike},
with a Courant algebroid structure \cite{R1ike}\cite{I4ike}\cite{HPike},
with a twisted Poisson structure \cite{KSike}\cite{I5ike},
with a Dirac structure \cite{KSSike},
with a generalized complex structure \cite{Zucike}\cite{I5ike}\cite{Pesike},
and so on.


\section{AKSZ formulation of Batalin-Vilkovisky formalism 
on Graded Bundles}
\label{sec:2}
\noindent
We explain general setting of the AKSZ formulation \cite{AKSZike} 
of the Batalin-Vilkovisky formalism for a general graded bundle
in the rather mathematical context.

Let $M$ be a smooth manifold in $d$ dimensions.
We consider a {\it supermanifold} $\Pi T^* M$.
Mathematically, $\Pi T^* M$, whose bosonic part is $M$, is defined as 
a cotangent bundle with reversed parity of the fiber.
That is, a base manifold $M$ has a Grassman even coordinate and 
the fiber of $\Pi T^* M$ has a Grassman odd coordinate.
We can introduce a grading.
The coordinate on the bass manifold has grade zero and 
the coordinate on the fiber has grade one.
This grading is called the {\it total degrees}.
Similarly, we can define $\Pi T M$ for a tangent bundle $T M$.

We can consider more general assignments for the degree of the fibers of 
$T^* M$ or $T M$.
For a nonnegative integer $p$, we define $T^* [p] M$, which is called 
a graded cotangent bundle.
$T^* [p] M$ is a cotangent bundle the degree of whose fiber is $p$.
If $p$ is odd, the fiber is Grassman odd, and if $p$ is even, 
the fiber is Grassman even. 
The coordinate on the bass manifold has the total degree zero and 
the coordinate on the fiber has the total degree $p$.
We define a graded tangent bundle $T [p] M$ in the same way.
%
For a general vector bundle $E$, a graded vector bundle $E[p]$ is defined 
in a similar way.
$E[p]$ is a vector bundle which degree of the fiber is shifted by $p$.
\footnote{Note that only the fiber is shifted and 
the base space is not shifted.}

We consider a Poisson manifold $N$ with a Poisson bracket $\{*,*\}$.
If we can construct a graded manifold $\tilde{N}$ from $N$,
 then the Poisson structure $\{*,*\}$ shifts to
a graded Poisson structure by grading of $\tilde{N}$.
The graded Poisson bracket is called an {\it antibracket} 
and denoted by $\sbvike{*}{*}$.
$\sbvike{*}{*}$ is graded symmetric and satisfies the graded Leibniz
rule and the graded Jacobi identity with respect to grading of the manifold. 
%
The antibracket $\sbvike{*}{*}$ with 
 the total degree $- n + 1$ satisfies the following identities:
\begin{eqnarray}
&& \sbvike{F}{G} = -(-1)^{(|F| + 1 - n)(|G| + 1 - n)} \sbvike{G}{F},
\nonumber \\
&& \sbvike{F}{G  H} = \sbvike{F}{G} H
+ (-1)^{(|F| + 1 - n)|G|} G \sbvike{F}{H},
\nonumber \\
&& \sbvike{F  G}{H} = F \sbvike{G}{H}
+ (-1)^{|G|(|H| + 1 - n)} \sbvike{F}{H} G,
\nonumber \\
&& (-1)^{(|F| + 1 - n)(|H| + 1 - n)} \sbvike{F}{\sbvike{G}{H}}
+ {\rm cyclic \ permutations} = 0,
\label{BVidentity}
\end{eqnarray}
where $F, G$ and $H$ are functions on $\tilde{N}$ and
$|F|, |G|$ and $|H|$ are the total degrees of functions respectively.
The graded Poisson structure is also called the {\it P-structure}.

Typical examples of a Poisson manifold $N$ are a cotangent bundle $T^* M$ 
and a vector bundle $E \oplus E^*$.
Two bundles have important roles in this paper.
An other example is a vector bundle $E$ with a Poisson structure on 
the fiber.
We consider these three bundles.

First we consider a cotangent bundle $T^* M$.
Since $T^* M$ has a natural symplectic structure,
we can define a Poisson bracket induced from the natural 
symplectic structure.
If we take a local coordinate $\phi^i$ on M and a local coordinate 
$B_i$ of the fiber, we can define a Poisson bracket as follows:
\footnote{We take Einstein's summation notation.} 
\begin{eqnarray}
\{ F, G \} \equiv 
F \frac{\rdike}{\partial \phi^i} 
\frac{\ldike }{\partial B_i}  G
- 
F \frac{\rdike }{\partial B_i} 
\frac{\ldike }{\partial \phi^i} G.
\label{ab0ike}
\end{eqnarray}
where $F$ and $G$ are a function on $T^* M$, and
${\rdike}/{\partial \varphi}$ and ${\ldike}/{\partial
\varphi}$ are the right and left differentiations
with respect to $\varphi$, respectively.
%
Next we shift the degree of fiber by $p$, 
and we consider the space $T^* [p]M$.
The Poisson structure changes to a graded Poisson structure.
The corresponding graded Poisson bracket is called the {\it antibracket}, 
$\sbvike{*}{*}$.
%
Let $\bphiike^i$ be a local coordinate of $M$
 and $\bb_{n-1,i}$ a basis of the fiber of $T^* [p] M$.
\footnote{We use bold notations for local coordinates of 
graded (super)vector bundles, while
we use nonbold notations for local coordinates of usual vector bundles.}
An {\it antibracket} $\sbvike{*}{*}$ 
on a cotangent bundle $T^* [p] M$ is represented as:
%
\begin{eqnarray}
\sbvike{F}{G} \equiv 
F \frac{\rdike}{\partial \bphiike^i} 
\frac{\ldike }{\partial \bb_{p,i}}  G
- 
F \frac{\rdike }{\partial \bb_{p,i}} 
\frac{\ldike }{\partial \bphiike^i} G.
\label{Poisson2ike}
\end{eqnarray}
The total degree of the $\sbvike{*}{*}$ is $-p$.

Next, we consider a vector bundle $E \oplus E^*$. 
There is a natural Poisson structure on the fiber of $E \oplus E^*$ induced 
from natural paring among $E$ and $E^*$.
If we take a local coordinate $A^a$ on the fiber of $E$ and 
$B_a$ on the fiber of $E^*$,
we define 
\begin{eqnarray}
\{ F, G \} \equiv 
F \frac{\rdike}{\partial A^a} 
\frac{\ldike }{\partial B_a}  G
- 
F \frac{\rdike }{\partial B_a} 
\frac{\ldike }{\partial A^a} G.
\label{abpike}
\end{eqnarray}
where $F$ and $G$ are a function on $E \oplus E^*$.
We shift the degrees of fibers of $E$ and $E^*$
to $E [p] \oplus E^*[q]$, where
$p$ and $q$ are positive integers.
The Poisson structure changes to a graded Poisson structure
$\sbvike{*}{*}$.
Let $\ba_p{}^{a}$ be a basis of the fiber of 
$E[p]$ and 
$\bb_{q,a}$ a basis of the fiber of 
$E^*[q]$.
An antibracket is represented as
\begin{eqnarray}
\sbvike{F}{G} \equiv 
F  \frac{\rdike}{\partial \ba_p{}^{a}} 
\frac{\ldike }{\partial \bb_{q,a}}  G
- (-1)^{p q}
F \frac{\rdike }{\partial \bb_{q,a}} 
\frac{\ldike }{\partial \ba_p{}^{a}}  G.
\label{Poisson4ike}
\end{eqnarray}
The total degree of the $\sbvike{*}{*}$ is $-p-q$.

Next, we consider a vector bundle $E$ with a Poisson structure on the fiber.
If we shift the degree of the fiber of $E$ to $E [p]$,
the Poisson structure changes to a graded Poisson structure
$\sbvike{*}{*}$.
Let $\ba_p{}^{a}$ be a basis of the fiber of 
$E[p]$,
An antibracket is represented as
\begin{eqnarray}
\sbvike{F}{G} \equiv 
F  \frac{\rdike}{\partial \ba_p{}^{a}} k^{ab} 
\frac{\ldike }{\partial \ba_p{}^{b}}  G,
\label{Poisson6ike}
\end{eqnarray}
where $F$ and $G$ are a function on $E[p]$ and 
$k^{ab}$ is a nondegenerate constant bivector 
induced from a (graded) Poisson structure.
The total degree of the antibracket $\sbvike{*}{*}$ is $-2p$.


Next we define a {\it Q-structure}.
A {\it Q-structure} is 
a function $S$ on a supermanifold $\tilde{N}$
which satisfies the classical master equation 
$\sbvike{S}{S} = 0$.
$S$ is called a {\it Batalin-Vilkovisky action}, or simply 
an {\it action}.
We require that $S$ satisfy the compatibility condition
\begin{eqnarray}
S \sbvike{F}{G} = \sbvike{S F}{G} + (-1)^{|F| +1} \sbvike{F}{SG},
\end{eqnarray}
where $F$ and $G$ are arbitrary functions and 
$|F|$ is the total degree of $F$.
$\sbvike{S}{F} = \brsike F$ generates an infinitesimal 
transformation(a Hamiltonian flow).
We call this a {\it BRST transformation}, which
coincides with the gauge transformation of the theory.

We define a (classical) {\it Batalin-Vilkovisky structure} as follows:
\begin{definition}
If a structure on a {\it supermanifold} 
has {\it P-structure} and {\it Q-structure}, 
it is called a {\it Batalin-Vilkovisky structure}.
\end{definition}


\section{Batalin-Vilkovisky Structures of Abelian Topological Sigma Models}
\subsection{BF case}
\noindent
In this section, 
we consider Batalin-Vilkovisky structures of topological sigma models.

Let $X$ be a base manifold in $n$ dimensions, with or without boundary, 
and $M$ be a target manifold in $d$ dimensions.
We denote $\phi$ a smooth map from $X$ to $M$.
%

We consider a {\it supermanifold} $\Pi T X$, whose bosonic part is $X$.
$\Pi T X$ is defined as 
a tangent bundle with reversed parity of the fiber.
We extend a smooth function $\phi$ to a function $\bphiike:\Pi T X \rightarrow 
M$.
$\bphiike$ is an element of $\Pi T^* X \otimes M$.
The {\it total degree} defined in the previous section is grading 
with respect to $M$.
We introduce a nonnegative integer grading on $\Pi T^* X$.
A coordinate on a bass manifold is zero and a coordinate on the fiber is one.
This grading is called the {\it form degrees}.
We denote ${\rm deg} F$ the form degree of a function $F$.
${\rm gh} F = |F| - {\rm \deg} F$ is called the {\it ghost number}.

First we consider a {\it P-structure} on $T^* [p] M$.
It is natural to take $p=n-1$ to construct
a Batalin-Vilkovisky structure in a topological sigma model.
In other words, 
the dimensions of $X$ labels the total degree 
of a Batalin-Vilkovisky structure on the supermanifold $T^* [p] M$.
We consider $T^* [n-1] M$ for an $n$-dimensional base manifold $X$.
Let $\bphiike^i$ be a local coordinate expression of $\Pi T^* X \otimes M$, 
where $i, j, k, \cdots$ are indices of local coordinate on $M$.
Let $\bb_{n-1,i}$ be a basis 
of sections of $\Pi T^* X \otimes \bphiike^*(T^* [n-1] M)$.
According to the discussion in the previous section, 
we can define an {\it antibracket} $\sbvike{*}{*}$ 
on a cotangent bundle $T^* [n-1] M$ as
%
\begin{eqnarray}
\sbvike{F}{G} \equiv 
F \frac{\rdike}{\partial \bphiike^i} 
\frac{\ldike }{\partial \bb_{n-1,i}}  G
- 
F \frac{\rdike }{\partial \bb_{n-1,i}} 
\frac{\ldike }{\partial \bphiike^i} G,
\label{cotangentPike}
\end{eqnarray}
where 
$F$ and $G$ are functions of $\bphiike^i$ and $\bb_{n-1,i}$.
We take a Darboux coordinate $\bphiike^i, \bb_{n-1,i}$, 
but it is for simplicity. 
We can take more general coordinates.
The total degree of the antibracket is $-n+1$.
If $F$ and $G$ are functionals of 
$\bphiike^i$ and $\bb_{n-1,i}$, 
we understand an antibracket is defined as 
\begin{eqnarray}
\sbvike{F}{G} \equiv 
\int_{X} 
F \frac{\rdike}{\partial \bphiike^i} 
\frac{\ldike }{\partial \bb_{n-1,i}}  G
- 
F \frac{\rdike }{\partial \bb_{n-1,i}} 
\frac{\ldike }{\partial \bphiike^i} G,
\end{eqnarray}
where the integration over $X$ is understood as that on the 
$n$-form part of the integrand. 
Through this article, 
we always understand an antibracket on two functionals in a similar 
manner and abbreviate this notation.

Next we consider a {\it P-structure} on $E \oplus E^*$.
Natural assignment of the total degrees
is nonnegative integers $p$ and $q$ are $p+q = n-1$.
That is, 
we consider $E[p] \oplus E^*[n-p-1]$, 
where $1 \leq p \leq n-2$.
We can naturally construct a topological sigma model in this case.
Let $\lfloor  x \rfloor$ be the floor function
which gives the largest integer less than or equal to $x$.
If $\lfloor \frac{n}{2}\rfloor \leq p \leq n-2$,
we identify 
$E[p] \oplus E^*[n-p-1]$ with the dual bundle 
$E^*[n-p-1] \oplus (E^*)^*[p]$.
Therefore 
$1 \leq p \leq \lfloor \frac{n-1}{2}\rfloor$ 
provides different structures of grading.

Let $\ba_p{}^{a_p}$ be a basis of sections of 
$\Pi T^* X \otimes \bphiike^*(E[p])$ and 
$\bb_{n-p-1,a_p}$ a basis of the fiber of 
$\Pi T^* X \otimes \bphiike^*(E^*[n-p-1])$.
From (\ref{Poisson4ike}), we can define an antibracket 
as
\begin{eqnarray}
\sbvike{F}{G} \equiv 
F  \frac{\rdike}{\partial \ba_p{}^{a_p}} 
\frac{\ldike }{\partial \bb_{n-p-1,a_p}}  G
- (-1)^{n p}
F \frac{\rdike }{\partial \bb_{n-p-1,a_p}} 
\frac{\ldike }{\partial \ba_p{}^{a_p}}  G.
\label{EEPike}
\end{eqnarray}


We want to consider various grading assignments for $E \oplus E^*$.
Because each assignment 
induces different Batalin-Vilkovisky structures.
In order to consider all independent assignments, 
we define the following bundle.
Let $E_p$ be $\lfloor \frac{n-1}{2}\rfloor$ 
series of vector bundles, 
where $1 \leq p \leq \lfloor \frac{n-1}{2}\rfloor$.
We consider $E_p[p] \oplus E_p^*[n-p-1]$ and consider a direct sum 
\begin{eqnarray}
\sum_{p=1}^{\lfloor \frac{n-1}{2} \rfloor}
E_p[p] \oplus E_p^*[n-p-1]
\label{totbundleike}
\end{eqnarray}
And we define a {\it P-structure}
on the graded vector bundle
\begin{eqnarray}
\left(\sum_{p=1}^{\lfloor \frac{n-1}{2} \rfloor}
E_p[p] \oplus E_p^*[n-p-1] \right) \oplus T^*[n-1] M.
\label{totspaceike}
\end{eqnarray}
A local (Darboux) coordinate expression of 
the antibracket $\sbvike{\cdot}{\cdot}$ is a sum of 
(\ref{cotangentPike}) and (\ref{EEPike}):
\begin{eqnarray}
\sbvike{F}{G} \equiv 
\sum_{p=0}^{\lfloor \frac{n-1}{2} \rfloor} 
F \frac{\rdike}{\partial \ba_p{}^{a_p}} 
\frac{\ldike }{\partial \bb_{n-p-1 \ a_p}} G
- (-1)^{n p}
F \frac{\rdike }{\partial \bb_{n-p-1 \ a_p}} 
\frac{\ldike }{\partial \ba_p{}^{a_p}} G.
\label{bfantibracketike}
\end{eqnarray}
where 
$p=0$ component is the antibracket (\ref{cotangentPike}) on 
the graded cotangent bundle $T^* [n-1] M$
and $\ba_0{}^{a_0} = \bphiike^i$.
Note that all terms of the antibracket have 
the total degree $-n+1$.
We can confirms that 
the antibracket (\ref{bfantibracketike}) satisfies the identity
(\ref{BVidentity}).


We construct a {\it Q-structure} on the bundle (\ref{totspaceike}).
A simplest and natural action for a topological sigma model is 
\begin{eqnarray}
S_0 = && \sum_{p=0}^{\lfloor \frac{n-1}{2} \rfloor} 
(-1)^{n-p} \int_{X} \bb_{n-p-1\ a_p} d \ba_p{}^{a_p},
\label{s0ike}
\end{eqnarray}
where $d$ is a exterior differential on $X$.
The total degree of $d$ is $1$ because 
the form degree is $1$ and we assign the ghost number $0$. 
The integration over $X$ is understood as that over the 
$n$-form part (the form degree $n$ part) of the integrand. 
This action is 
analogy of a fundamental form $\theta = p_i d q^i$ for a
symplectic form $\omega$, which has $\omega = - d \theta$,
therefore this action is directly 
derived from the {\it P-structure} on the graded bundle.
The total degree of $S_0$ is $n$.
It is called an {\it abelian BF theory in $n$ dimensions}.

$S_0$ defines a {\it Q-structure}, since 
we can easily confirm that $S_0$ 
satisfies the classical master equation:
\begin{eqnarray}
(S_0, S_0) = 0,
\end{eqnarray}

'Abelian' means that the theory has a $U(1)$ gauge symmetry.
The BRST transformation (the gauge symmetry) is defined as
\begin{eqnarray}
\brsike_0 \bPhiike \equiv (S_0, \bPhiike) = d \bPhiike,
\label{brstphiike}
\end{eqnarray}
where $\bPhiike$ is an arbitrary section of a total bundle
$
\left(\sum_{p=1}^{\lfloor \frac{n-1}{2} \rfloor}
E_p[p] \oplus E_p^*[n-p-1] \right) \oplus T^*[n-1] M
$.
$\brsike_0^2 = 0$ is satisfied from $(S_0, S_0) = 0$, 
which is consistent to $d^2 = 0$.

\subsection{Chern-Simons with BF case}
\noindent
For a vector bundle $E$ with a Poisson structure on the fiber,
we can construct an another topological sigma model
if $n$ is odd. 
%
%
We consider a graded vector bundle $E[q]$ for $E$, 
which degree of the fiber is shifted by $q$.
Let $\ba_q{}^{a_q}$ be a basis of sections of 
$\Pi T^* X \otimes \bphiike^*(E[q])$.
From the equation (\ref{Poisson6ike}), 
we can define an antibracket as 
\begin{eqnarray}
\sbvike{F}{G} \equiv 
F  \frac{\rdike}{\partial \ba_q{}^{a}} k^{ab} 
\frac{\ldike }{\partial \ba_q{}^{b}}  G,
\label{EPike}
\end{eqnarray}

If we take a nonnegative integer $q = \frac{n-1}{2}$,
then the total degrees of 
(\ref{cotangentPike}), (\ref{EEPike}) and (\ref{EPike})
are all $-n+1$.
Thus we can consider a combined {\it P-structure}.

We consider a direct sum of $E[q] = E\left[\frac{n-1}{2}\right]$
with $E_p[p] \oplus E_p^*[n-p-1]$
\begin{eqnarray}
\left(\sum_{p=1}^{\frac{n-3}{2}}
E_p[p] \oplus E_p^*[n-p-1] \right) \oplus E\left[\frac{n-1}{2}\right],
\label{totbcsike}
\end{eqnarray}
where we have absorbed 
$E_p[p] \oplus E_p^*[n-p-1]$
to $E\left[\frac{n-1}{2}\right]$ for $p = \frac{n-1}{2}$.
We can define a natural {\it P-structure}
on the graded vector bundle
\hfil\break
$
\left(\sum_{p=1}^{\frac{n-3}{2}}
E_p[p] \oplus E_p^*[n-p-1] \right) \oplus E\left[\frac{n-1}{2}\right]
\oplus T^*[n-1] M
$.
%
A local (Darboux) coordinate expression of 
the antibracket $\sbvike{\cdot}{\cdot}$ is a sum of 
(\ref{cotangentPike}), (\ref{EEPike}) and 
(\ref{EPike}):
\begin{eqnarray}
\sbvike{F}{G} \equiv &&
\sum_{p=0}^{\frac{n-3}{2}} 
F \frac{\rdike}{\partial \ba_p{}^{a_p}} 
\frac{\ldike }{\partial \bb_{n-p-1 \ a_p}} G
- (-1)^{n p}
F \frac{\rdike }{\partial \bb_{n-p-1 \ a_p}} 
\frac{\ldike }{\partial \ba_p{}^{a_p}} G
\nonumber \\
&&
+ F  \frac{\rdike}{\partial \ba_q{}^{a}} k^{ab} 
\frac{\ldike }{\partial \ba_q{}^{b}}  G,
\label{csantiike}
\end{eqnarray}
where 
$p=0$ component is the antibracket on 
the graded cotangent bundle $T^* [n-1] M$
and $\ba_0{}^{a_0} = \bphiike^i$.
$q = \frac{n -1}{2}$ is needed for 
all the terms to have the total degree $-n+1$.


A simplest and natural action $S_0$ is 
\begin{eqnarray}
S_0 = && \sum_{p=0}^{\frac{n-3}{2}} 
(-1)^{n-p} \int_{X} \bb_{n-p-1\ a_p} d \ba_p{}^{a_p}
+ \int_{X} \frac{k^{ab}}{2} 
\ba_{\frac{n-1}{2}\ a} d \ba_{\frac{n-1}{2}\ b},
\label{s0csike}
\end{eqnarray}
%
%
The second term 
is 
called an {\it abelian Chern-Simons theory} in $n$ dimensions.
$S_0$ defines a {\it Q-structure}, since 
we can easily confirm that 
$S_0$ satisfies the classical master equation:
\begin{eqnarray}
(S_0, S_0) = 0,
\end{eqnarray}
This action also has a $U(1)$ gauge symmetry. 
The BRST transformation is
\begin{eqnarray}
\brsike_0 \bPhiike = (S_0, \bPhiike) = d \bPhiike,
\label{brstcdphiike}
\end{eqnarray}
where $\bPhiike$ is an arbitrary section of the total bundle.


\section{Deformation}
\noindent
In this section, we consider deformation of Batalin-Vilkovisky structures.
Deformation means deformation of the {\it Q-structure} 
for a fixed {\it P-structure}.
We consider local deformation from the settled point (\ref{s0ike})
or (\ref{s0csike}) of the moduli space.

First we consider the BF case of the bundle
$
\left(\sum_{p=0}^{\lfloor \frac{n-1}{2} \rfloor}
E_p[p] \oplus E_p^*[n-p-1] \right) \oplus T^*[n-1] M
$.
The beginning {\it Q-structure} is the equation (\ref{s0ike}). 
We deform this Batalin-Vilkovisky action $S_0$ to 
\begin{eqnarray}
S = S_0 + g S_1,
\label{lapSike}
\end{eqnarray}
under 
the condition that $S$ also satisfies the classical master equation:
\begin{eqnarray}
\sbvike{S}{S} = 0,
\label{smasterike}
\end{eqnarray}
where $g$ is a deformation parameter and $S_1$ represents
all the deformation terms, which are functional on $X$ and 
functions on 
$
\left(\sum_{p=0}^{\lfloor \frac{n-1}{2} \rfloor}
E_p[p] \oplus E_p^*[n-p-1] \right) \oplus T^*[n-1] M
$.
We require that $S$ is the total degree $n$. It is equivalent 
that the ghost number is zero ${\rm gh} S = |S| - {\rm deg} S = 0$.
This condition is physically necessary, 
though we can relax this condition mathematically.
If two deformations $S$ and $S^{\prime}$ satisfy 
$S^{\prime} = S + \sbvike{S_0}{T} = S + \brsike_0{T}$ 
for some functional $T$, 
two Batalin-Vilkovisky structures are equivalent.
Therefore the problem is to look for the total degree 
$n$ cohomology class 
\hfil\break
$H^n \left(\Pi T^* X \otimes
\left(\sum_{p=0}^{\lfloor \frac{n-1}{2} \rfloor}
E_p[p] \oplus E_p^*[n-p-1] \right) \oplus T^*[n-1] M \right)$.

Substituting $S = S_0 + g S_1$ to (\ref{smasterike}), 
we obtain $g$ expansion 
\begin{eqnarray}
\sbvike{S_0}{S_0} + 2 g \sbvike{S_0}{S_1}
+ g^2 \sbvike{S_1}{S_1} = 0,
\label{expabike}
\end{eqnarray}

Zero-th order of the equation (\ref{expabike}), 
$\sbvike{S_0}{S_0} = 0$, is already satisfied. 

First order 
is $\sbvike{S_0}{S_1} = 0$.
Because of the equation (\ref{brstphiike}),
$S_1$ is the integration of 
an arbitrary function $F$ of all fundamental superfields $\bPhiike$
which are $\ba_{p}{}^{a_{p}}$ or $\bb_{q b_{q}}$,
and their derivatives $d \bPhiike$:
\begin{eqnarray}
S_1 = && \int_{X} F(\bPhiike, d \bPhiike).
\end{eqnarray}
For simplicity, we assume that there is no boundary contribution to $S$,
that is, integration of total derivative terms on $X$ is always zero.
This corresponds to assume that there is no obstruction of deformation.
Then we can prove the following theorem:
\begin{theorem}
Assume there is no boundary contribution on $X$, i.e.
$\int_{X} d G(\bPhiike) = 0$ for any function $G$.
If a monomial of $F(\bPhiike, d \bPhiike)$ includes at least one $d \bPhiike$, 
$\int_{X} F(\bPhiike, d \bPhiike)$ is $\brsike_0$-exact.
\end{theorem}
\begin{proof}
We can assume that
$\int_{X} F(\bPhiike, d \bPhiike) 
= \sum_{p=0}^{n-1} \int_{X} F_{n-p-1} d G_{p}$,
where 
$F_{n-p-1}$ are functions with the form degree 
$n-p-1$ and $G_{p}$ are functions with the form degree 
$p$.
From (\ref{brstphiike}), 
we obtain
\begin{eqnarray}
&& \delta_0 F_{0} = 0, 
\nonumber \\
&& \delta_0 F_{n-p-1} = d F_{n-p-2} \quad \mbox{for $-1 \leq p \leq n-2$},
\nonumber \\
&& d F_{n} = 0, 
\nonumber \\
&& \delta_0 G_{0} = 0, 
\nonumber \\
&& \delta_0 G_{p} = d G_{p-1} \quad \mbox{for $1 \leq p \leq n$},
\nonumber \\
&& d G_{n} = 0, 
\label{expbrstike}
\end{eqnarray}
For even $p$, adjoining two terms are combined as
\begin{eqnarray}
&& F_{n-p-1} d G_{p} + F_{n-p-2} d G_{p+1}
\nonumber \\
&& =
(-1)^{n-p-1} \delta_0( F_{n-p-1} G_{p+1} )
- (-1)^{n-p-1} d( F_{n-p-2} G_{p+1} )
\label{2termsike}
\end{eqnarray}
using the relations (\ref{expbrstike}). 
Thus integration of these two terms is $\delta_0$-exact.

If $n$ is even, 
$S_1 = \sum_{p=0}^{n-1} \int_{X} F_{n-p-1} d G_{p}$ has even numbers of terms,
therefore 
we combine each two term like (\ref{2termsike}) and
we can confirm that $S_1$ is $\delta_0$-exact.

If $n$ is odd, the last term 
$F_{0} d G_{n-1}$ remains.
However this term is $\delta_0$-exact because
$F_{0} d G_{n-1} = \delta_0 (F_{0} G_{n})$.
Therefore $S_1$ is $\delta_0$-exact.
\smartqed\qed 
\end{proof}
Therefore since the nontrivial $S_1$ cohomological class
does not include $d$, we can take 
$S_1 = \int_{X} F(\bPhiike)$.
Concretely we can express
\begin{eqnarray}
S_1 = && \sum_{p(1), \cdots, p(k), q(1), \cdots, q(l)} \int_{X} 
F_{p(1) \cdots p(k), q(1) \cdots q(l) \ a_{p(1)} \cdots a_{p(k)}}
{}^{b_{q(1)} \cdots b_{q(l)}}(\ba_0{}^{a_0})
\nonumber \\
&& \times \ba_{p_1}{}^{a_{p(1)}} \cdots \ba_{p_k}{}^{a_{p(k)}}
\bb_{q_1 b_{q(1)}} \cdots \bb_{q_l b_{q(l)}},
\label{s1ike}
\end{eqnarray}
where $F_{p(1) \cdots p(k), q(1) \cdots q(l) \ a_{p(1)} \cdots a_{p(k)}}
{}^{b_{q(1)} \cdots b_{q(l)}}(\ba_0{}^{a_0})$ 
is a function of $\ba_0{}^{a_0}$
and $p(r) \neq 0, q(s) \neq 0$ for $r = 1, \cdots, k, s = 1, \cdots, l$.
From the ${\rm deg} S = n$ condition, we obtain
$\sum_{\alpha=0}^k
 |\ba_{p_\alpha}{}^{a_{p(\alpha)}}| +
\sum_{\beta=0}^l |\bb_{q_\beta b_{q(\beta)}}| = n$.

Second order of the equation (\ref{expabike})
\begin{eqnarray}
\sbvike{S_1}{S_1} = 0,
\label{s1s1ike}
\end{eqnarray}
imposes conditions on functions
$F_{p(1) \cdots p(k), q(1) \cdots q(l) \ a_{p(1)} \cdots a_{p(k)}}
{}^{b_{q(1)} \cdots b_{q(l)}}(\ba_0{}^{a_0})$.
These conditions determine the mathematical structure 
of a Batalin-Vilkovisky structure.

We call a resulting field theory $S=S_0+g S_1$ a {\it nonlinear gauge theory}
in $n$ dimensions.

Next we consider the Chern-Simons with BF case 
\hfil\break
$
\left(\sum_{p=1}^{\frac{n-3}{2}}
E_p[p] \oplus E_p^*[n-p-1] \right) \oplus E\left[\frac{n-1}{2}\right].
$
We make a similar discussion with the BF case.
%
We consider deformation of the Batalin-Vilkovisky action (\ref{s0csike}) to 
(\ref{lapSike}), $S = S_0 + g S_1$,
under 
the condition that $S$ also satisfies the classical master equation
(\ref{smasterike}).
We obtain $S_1$ by using Theorem 1.
$S_1$ has a similar expression with
(\ref{s1ike}) but has different field contents.
The second order of the equation (\ref{expabike})
\begin{eqnarray}
\sbvike{S_1}{S_1} = 0,
\end{eqnarray}
imposes conditions on functions
$F_{p(1) \cdots p(k), q(1) \cdots q(l) \ a_{p(1)} \cdots a_{p(k)}}
{}^{b_{q(1)} \cdots b_{q(l)}}(\ba_0{}^{a_0})$.
%


\section{Deformation in lower dimensions}
\noindent
In this section, we concretely 
analyze the algebraic and geometric structure of 
deformation of a topological sigma model in lower dimensional $X$.

We can easily find that we cannot obtain nontrivial deformation in case of 
the total degree $n=1$, thus we cannot obtain a nontrivial structure.
\subsection{$n=2$}
\noindent
We analyze the algebraic structure of the total degree $n=2$ 
topological sigma model.
In two dimensions, the total graded bundle (\ref{totspaceike}) is
$
T^*[1] M = \Pi T^* M 
$.

(\ref{lapSike}) 
under (\ref{s1ike}) is 
\begin{eqnarray}
&& S = S_0 +g S_1,  \nonumber \\
&& S_0 = \int_{X}
\bb_{1 i} d \bphiike^i, 
\qquad
S_1 = \int_{\Sigma} \frac{1}{2} f^{ij}(\bphiike)
\bb_{1 i} \bb_{1 j},
\label{2dactionike}
\end{eqnarray}
where $i, j, \cdots$ are indices of 
a local coordinate expressions on
$
T^*[1] M
$.
and 
we rewrite the notations as 
$\bphiike^i = \ba_0{}^i$ and 
$\frac{1}{2}f^{ij}(\bphiike) = F_{,11}{}^{ij}(\ba_0{})$.
This topological sigma model 
is known as {\it the Poisson sigma model} \cite{IIike}\cite{SSike}.
If we substitute this $S_1$ to the condition (\ref{s1s1ike}),
we obtain the geometric structure of the action.
We obtain the following identity on $f^{ij}$: 
\begin{eqnarray}
f^{kl} \frac{\ldike}{\partial \bphiike^l} f^{ij}
+ f^{il} \frac{\ldike}{\partial \bphiike^l} f^{jk}
+ f^{jl} \frac{\ldike}{\partial \bphiike^l} f^{ki}= 0.
\label{2dJacobiike}
\end{eqnarray}
If we restrict the identity on $X$ (i.e. a ghost number zero sector), 
(\ref{2dJacobiike}) reduces to
\begin{eqnarray}
f^{kl}(\phi) \frac{\partial f^{ij}(\phi)}{\partial \phi^l}
+ f^{il}(\phi) \frac{\partial f^{jk}(\phi)}{\partial \phi^l}
+ f^{jl}(\phi) \frac{\partial f^{ki}(\phi)}{\partial \phi^l}= 0.
\label{2dJacobi2ike}
\end{eqnarray}
%
Under the identity (\ref{2dJacobiike}), 
$- f^{ij}$ defines a Poisson
structure as
\begin{eqnarray}
\{F(\phi),G(\phi) \} \equiv - f^{ij}(\phi)
\frac{\partial F}{\partial \phi^i}\frac{\partial G}{\partial \phi^j},
\label{pbracketike}
\end{eqnarray}
on the space $M$.
Conversely, if we consider the Poisson structure $- f^{ij}$ on $M$, which
satisfies the identity (\ref{2dJacobiike}), we can define the action 
(\ref{2dactionike}) consistently.
This Poisson bracket (\ref{pbracketike}) is directly constructed 
if we restrict the derived bracket \cite{Kosike} 
of a Batalin-Vilkovisky structure:
\begin{eqnarray}
\{F(\bphiike),G(\bphiike) \} = - 
F \frac{\rdike}{\partial \bphiike^i}f^{ij}(\bphiike)\frac{\ldike}{\partial \bphiike^j} G
= \sbvike{\sbvike{S}{F}}{G},
\end{eqnarray}
to $M$ and $X$, i.e. to the ghost number zero sector.

%

The Batalin-Vilkovisky structure of this theory has a structure of 
the {\it Lie algebroid} \
$T^* [1] \calMike$ \cite{LOike}\cite{Oike},
where $\calMike$ is a space of a (smooth) map from $\Pi T X$
to a target space $M$.

A {\it Lie algebroid} is a generalization of a bundle 
of a Lie algebra over a base manifold $\calMike$.
A Lie algebroid over a manifold is a vector bundle 
$\calEike \rightarrow \calMike$ with a Lie algebra structure on the 
space of the sections $\Gamma(\calEike)$ defined by the 
Lie bracket $[e_1, e_2], \quad e_1, e_2 \in \Gamma(\calEike)$
and a bundle map (the anchor)
$\rho: \calEike \rightarrow T\calMike$ satisfying the following properties:
\begin{eqnarray}
&& 1, \quad \rm{For \ any} \quad e_1, e_2 \in \Gamma(\calEike), \quad
[\rho(e_1), \rho(e_2)] = \rho([e_1, e_2]),
\nonumber \\
&& 2, \quad \rm{For \ any} \quad e_1, e_2 \in \Gamma(\calEike), 
\ F \in C^{\infty}(\calMike), 
\nonumber \\ 
&& \qquad [e_1, F e_2] = F [e_1, e_2] + (\rho(e_1) F) e_2,
  \label{liealgdef}
\end{eqnarray}

We can derive the bracket of a Lie algebroid from the antibracket
and the Batalin-Vilkovisky structure of the $n=2$ theory.
In our case, 
a vector bundle $\calEike$ is a cotangent bundle $T^*[1]\calMike$.
%
The Lie bracket of the two sections $e_1$ and $e_2$ is defined as 
a derived bracket of the antibracket by
\begin{eqnarray}
[e_1, e_2] \equiv \sbvike{\sbvike{S}{e_1}}{e_2},
\end{eqnarray}
and the anchor is defined by
\begin{eqnarray}
\rho(e) F(\bphiike) \equiv \sbvike{e}{\sbvike{S}{F(\bphiike)}}.
\end{eqnarray}
Then we can confirm $[e_1, e_2] = - [e_2, e_1]$,  
from $\sbvike{e_1}{e_2}=0$ and the graded Jacobi identity of the antibracket.
Similarly, a Lie algebroid conditions 1 and 2 on the bracket 
$[\cdot, \cdot]$ and the anchor map $\rho$ 
are obtained from the classical master equation $\sbvike{S}{S}=0$.
From this derived bracket,
we directly obtain 
the following ``noncommutative'' relation on the coordinates:
\begin{eqnarray}
[\bphiike^i, \bphiike^j] = - f^{ij}(\bphiike).
\end{eqnarray}
The anchor is a differentiation on functions of $\bphiike$ as
\begin{eqnarray}
\rho(\bphiike^i) F(\bphiike) = - f^{ij}(\bphiike) 
\frac{\ldike}{\partial \bphiike^j} F(\bphiike).
\end{eqnarray}

\subsection{$n=3$}
\subsubsection{BF case}
\noindent
We analyze a nonlinear gauge theory in three dimensions.
In this case, the theory defines the topological open 2-brane as a
sigma model \cite{I2ike}\cite{Pike}.

We consider a supermanifold $\Pi T X$ whose ghost number zero part is 
a three-dimensional manifold $X$.
A base space $\calMike$ is the space of a (smooth) map $\bphiike$ 
from $\Pi T X$ to a target space $M$.
In $n=3$, the total bundle (\ref{totbundleike}) is 
a vector bundle $E[1] \oplus E^*[1]$.
From (\ref{totspaceike}), we consider 
$(E[1] \oplus E^*[1]) \oplus T^*[2]M$.
We obtain an antibracket $\sbvike{\cdot}{\cdot}$ 
on the space of sections of 
$\Pi T^* X \otimes \bphiike^*(E[1] \oplus E^*[1] \oplus T^*[2]M)$
if we set $n=3$ in (\ref{bfantibracketike}).

In order to write down the Batalin-Vilkovisky action $S$, 
we take a local basis on 
$\Pi T^* X \otimes \bphiike^*(E[1] \oplus E^*[1])$
as $\ba_{1}{}^a, \bb_{1\ a}$, 
which are Darboux coordinates
such that $\sbvike{\ba_1{}^a}{\ba_1{}^b} = 
\sbvike{\bb_{1a}}{\bb_{1b}} = 0$ and
$\sbvike{\ba_1{}^a}{\bb_{1b}} = \delta^{a}{}_b$.
Moreover we introduce $\bb_{2 i}$
a section of $\Pi T^* X \oplus \bphiike^*(T^*[2]M)$.

The total action (\ref{lapSike}) 
under (\ref{s1ike}) is 
\begin{eqnarray}
&& S = S_0 +g S_1,  \nonumber \\
&& S_0 = \int_{X}
[- \bb_{2 i} d \bphiike{}^i + \bb_{1 a} d \ba_1{}^a], 
\nonumber \\
&& S_1 = \int_{X} 
[f_1{}_a{}^i(\bphiike) \ba_1{}^a \bb_{2 i} 
+ f_2^{ib}(\bphiike) \bb_{2 i} \bb_{1 b}
+ \frac{1}{3!} f_{3abc}(\bphiike) \ba_1{}^a \ba_1{}^b \ba_1{}^c
\nonumber \\
&& 
+ \frac{1}{2} f_{4ab}{}^c(\bphiike) \ba_1{}^a \ba_1{}^b \bb_{1 c}
+ \frac{1}{2} f_{5a}{}^{bc}(\bphiike) \ba_1{}^a \bb_{1 b} \bb_{1 c}
\nonumber \\
&& 
+ \frac{1}{3!} f_6{}^{abc}(\bphiike) \bb_{1 a} \bb_{1 b} \bb_{1 c}],
\label{3daction}
\end{eqnarray}
where we set
$f_1{}_a{}^i = F_{1,2 a}{}^i$,
$f_2^{ib} = F_{,21}{}^{ib}$,
$\frac{1}{3!}f_{3abc} = F_{111, abc}$,
$\frac{1}{2}f_{4ab}{}^c = F_{11,1 ab}{}^c$,
$\frac{1}{2}f_{5a}{}^{bc} = F_{1,11 a}{}^{bc}$,
$\frac{1}{3!}f_6{}^{abc} = F_{,111 }{}^{abc}$,
for clarity. 
The condition of the classical
master equation (\ref{s1s1ike}) imposes the following identities on six
$f_i$'s, $i=1, \cdots, 6$\cite{I2ike}:
\begin{eqnarray}
&& 
f_{1}{}_e{}^i f_{2}{}^{je} + f_{2}{}^{ie} f_{1}{}_e{}^j =
0,
\nonumber \\
&& 
- \left( \frac{\ldike}{\partial \bphiike^j} f_{1}{}_c{}^i 
\right) f_{1}{}_b{}^j 
+ \left( \frac{\ldike}{\partial \bphiike^j} f_{1}{}_b{}^i
\right) f_1{}_c {}^j
+ f_1{}_e{}^i f_{4bc}{}^e + f_{2}{}^{ie} f_{3ebc} = 0, 
\nonumber \\
&&
 f_1{}_b{}^j \left( \frac{\ldike}{\partial \bphiike^j} f_{2}{}^{ic} \right)
- f_{2}{}^{jc} \left( \frac{\ldike}{\partial \bphiike^j}  f_1{}_b{}^i \right)
+ f_1{}_e{}^i f_{5b}{}^{ec} - f_{2}{}^{ie} f_{4eb}{}^c = 0, 
\nonumber \\
&&
- f_{2}{}^{jb} \left( \frac{\ldike}{\partial \bphiike^j} f_{2}{}^{ic} \right)
+ f_{2}{}^{jc} \left( \frac{\ldike}{\partial \bphiike^j} f_{2}{}^{ib} \right)
+ f_1{}_e{}^i f_{6}^{ebc} + f_{2}{}^{ie} f_{5e}{}^{bc} = 0, 
\nonumber \\
&&
- f_1{}_{[a}{}^j \left( \frac{\ldike}{\partial \bphiike^j} f_{4bc]}{}^d \right)
+ f_{2}{}^{jd} \left( \frac{\ldike}{\partial \bphiike^j} f_{3abc} \right)
+ f_{4e[a}{}^d f_{4bc]}{}^{e} + f_{3e[ab} f_{5c]}{}^{de} =
0,
\nonumber \\
&&
- f_1{}_{[a}{}^j \left( \frac{\ldike}{\partial \bphiike^j} f_{5b]}{}^{cd} \right)
- f_{2}{}^{j[c} \left( \frac{\ldike}{\partial \bphiike^j} f_{4ab}{}^{d]} \right)
\nonumber \\
&& 
+ f_{3eab} f_6{}^{ecd} 
+ f_{4e[a}{}^{[d} f_{5b]}{}^{c]e} 
+ f_{4ab}{}^e f_{5e}{}^{cd} = 0, 
\nonumber \\
&&
- f_1{}_a{}^j \left( \frac{\ldike}{\partial \bphiike^j} f_{6}{}^{bcd} \right)
+ f_{2}{}^{j[b} \left( \frac{\ldike}{\partial \bphiike^j} f_{5a}{}^{cd]} \right)
+ f_{4ea}{}^{[b} f_6{}^{cd]e} + f_{5e}{}^{[bc} f_{5a}{}^{d]e} = 0, 
\nonumber \\
&&
- f_{2}{}^{j[a} \left( \frac{\ldike}{\partial \bphiike^j} f_{6}{}^{bcd]} \right)
+ f_6{}^{e[ab} f_{5e}{}^{cd]} = 0, 
\nonumber \\
&&
- f_1{}_{[a}{}^j \left( \frac{\ldike}{\partial \bphiike^j} f_{3bcd]} \right)
+ f_{4[ab}{}^{e} f_{3cd]e} = 0,
\label{3dJacobiike}
\end{eqnarray}
where $[\cdots]$ on the indices represents the antisymmetrization of
them, e.g., $\Phi_{[ab]} = \Phi_{ab} - \Phi_{ba}$.

If we restrict fields to $X$, 
(\ref{3dJacobiike}) reduces to the following identities:
\begin{eqnarray}
&& 
f_{1}{}_e{}^i f_{2}{}^{je} + f_{2}{}^{ie} f_{1}{}_e{}^j = 0, 
\nonumber \\
&& 
- \frac{\partial f_{1}{}_c{}^i}{\partial \phi^j} f_{1}{}_b{}^j
+ \frac{\partial f_{1}{}_b{}^i}{\partial \phi^j} f_1{}_c {}^j
+ f_1{}_e{}^i f_{4bc}{}^e + f_{2}{}^{ie} f_{3ebc} = 0, 
\nonumber \\
&&
 f_1{}_b{}^j \frac{\partial f_{2}{}^{ic}}{\partial \phi^j} 
- f_{2}{}^{jc} \frac{\partial f_1{}_b{}^i}{\partial \phi^j} 
+ f_1{}_e{}^i  f_{5b}{}^{ec} - f_{2}{}^{ie} f_{4eb}{}^c = 0, 
\nonumber \\
&&
- f_{2}{}^{jb} \frac{\partial f_{2}{}^{ic}}{\partial \phi^j} 
+ f_{2}{}^{jc} \frac{\partial f_{2}{}^{ib}}{\partial \phi^j} 
+ f_1{}_e{}^i f_{6}^{ebc} + f_{2}{}^{ie} f_{5e}{}^{bc} = 0, 
\nonumber \\
&&
- f_1{}_{[a}{}^j \frac{\partial f_{4bc]}{}^d}{\partial \phi^j} 
+ f_{2}{}^{jd} \frac{\partial f_{3abc}}{\partial \phi^j} 
+ f_{4e[a}{}^d f_{4bc]}{}^{e} + f_{3e[ab} f_{5c]}{}^{de} = 0, 
\nonumber \\
&&
- f_1{}_{[a}{}^j \frac{\partial f_{5b]}{}^{cd}}{\partial \phi^j} 
- f_{2}{}^{j[c} \frac{\partial f_{4ab}{}^{d]}}{\partial \phi^j} 
+ f_{3eab} f_6{}^{ecd} 
+ f_{4e[a}{}^{[d} f_{5b]}{}^{c]e} + f_{4ab}{}^e f_{5e}{}^{cd} = 0, 
\nonumber \\
&&
- f_1{}_a{}^j \frac{\partial f_{6}{}^{bcd}}{\partial \phi^j} 
+ f_{2}{}^{j[b} \frac{\partial f_{5a}{}^{cd]}}{\partial \phi^j} 
+ f_{4ea}{}^{[b} f_6{}^{cd]e} + f_{5e}{}^{[bc} f_{5a}{}^{d]e} = 0, 
\nonumber \\
&&
- f_{2}{}^{j[a} \frac{\partial f_{6}{}^{bcd]}}{\partial \phi^j} 
+ f_6{}^{e[ab} f_{5e}{}^{cd]} = 0, 
\nonumber \\
&&
- f_1{}_{[a}{}^j \frac{\partial f_{3bcd]}}{\partial \phi^j} 
+ f_{4[ab}{}^{e} f_{3cd]e} = 0.
\label{3dJacobi2ike}
\end{eqnarray}


The algebraic structure (\ref{3dJacobiike}) (or (\ref{3dJacobi2ike}))
is {\it a Courant algebroid}.
{\it A Courant algebroid} is introduced by Courant 
in order to analyze the Dirac structure as a generalization of a
Lie algebra of vector fields
on a vector bundle \cite{Cike}\cite{LWXike}.
The Batalin-Vilkovisky structure on a Courant algebroid 
is first analyzed in \cite{R1ike}.

{\it A Courant algebroid} is a vector bundle $\calEike \rightarrow \calMike$
and has a nondegenerate symmetric bilinear form
$\bracketike{\cdot}{\cdot}$ 
on the bundle, a bilinear operation $\circ$ on $\Gamma(\calEike)$ (the
space of sections on $\calEike$), and a bundle map (called the anchor) 
$\rho: \calEike \rightarrow T\calMike$ satisfying the following properties:
%
\begin{eqnarray}
&& 1, \quad e_1 \circ (e_2 \circ e_3) = (e_1 \circ e_2) \circ e_3 
+ e_2 \circ (e_1 \circ e_3), \nonumber \\
&& 2, \quad \rho(e_1 \circ e_2) = [\rho(e_1), \rho(e_2)], \nonumber \\
&& 3, \quad e_1 \circ F e_2 = F (e_1 \circ e_2)
+ (\rho(e_1)F)e_2, \nonumber \\
&& 4, \quad e_1 \circ e_2 = \frac{1}{2} {\cal D} \bracketike{e_1}{e_2},
\nonumber \\ 
&& 5, \quad \rho(e_1) \bracketike{e_2}{e_3}
= \bracketike{e_1 \circ e_2}{e_3} + \bracketike{e_2}{e_1 \circ e_3},
   \label{courantdefike}
\end{eqnarray}
where 
$e_1, e_2$ and $e_3$ are sections of $\calEike$, and $F$ is a function on 
$\calMike$;
${\cal D}$ is a map from functions on $\calMike$ to $\Gamma(\calEike)$ and is
defined by
$\bracketike{{\cal D}F}{e} = \rho(e) F$.
%

In our nonlinear gauge theory, 
$\calMike$ is the space of a map $\bphiike$ from $\Pi T X$ to $M$
and $\calEike$ is the space of sections of 
$\Pi T^* X \oplus \bphiike^*(E[1] \oplus E^*[1])$.
Let $e^a$ be a local basis of sections of $\calEike$.
We take $e^a = \ba_{1}{}^a$ or $\bb_{1\ a}$, 
which are a $E[1]$ component or a $E^*[1]$ component respectively.
We define a graded symmetric bilinear form 
$\bracketike{\cdot}{\cdot}$, a bilinear operation $\circ$, a
bundle map $\rho$ and ${\cal D}$ from the antibracket as follows:
\begin{eqnarray}
&& e^a \circ e^b \equiv \sbvike{\sbvike{S}{e^a}}{e^b}, \nonumber \\
&& \bracketike{e^a}{e^b} \equiv \sbvike{e^a}{e^b}, \nonumber \\
&& \rho(e^a) F(\bphiike) \equiv \sbvike{e^a}{\sbvike{S}{F(\bphiike)}}, \nonumber \\
&& {\cal D}(*) \equiv \sbvike{S}{*}.
  \label{corresbase2ike}
\end{eqnarray}
Then we can confirm that the classical master equation 
$\sbvike{S}{S} =0$, which derive 
the identity (\ref{3dJacobiike}) on structure functions $f$'s, 
is equivalent to the conditions $1$ to $5$ of 
the equation (\ref{courantdefike}).
\cite{R1ike}\cite{I4ike}
We calculate the operations 
$\circ$ and $\rho$ on the basis as follows:
\begin{eqnarray}
&& \ba_1{}^a \circ \ba_1{}^b = - f_{5c}{}^{ab}(\bphiike) \ba_1{}^c
- f_{6}{}^{abc}(\bphiike) \bb_{1c}, \nonumber \\
&& \ba_1{}^a \circ \bb_{1b} = - f_{4bc}{}^{a}(\bphiike) \ba_1{}^c
+ f_{5b}{}^{ac}(\bphiike) \bb_{1c}, \nonumber \\
&& \bb_{1a} \circ \bb_{1b} = - f_{3abc}(\bphiike) \ba_1{}^c
- f_{4ab}{}^{c}(\bphiike) \bb_{1c}, \nonumber \\
&& \rho(\ba_1{}^a) \bphiike^i = - f_{2}{}^{ia}(\bphiike),
\nonumber \\
&& \rho(\bb_{1a}) \bphiike^i = - f_{1a}{}^i(\bphiike).
  \label{bfcourant}
\end{eqnarray}
%
%
%
This topological sigma model defines a Courant algebroid structure on 
the space $E[1] \oplus E^*[1]$.
We call this model as {\it the Courant sigma model}.

\subsubsection{Chern-Simons with BF case}
\noindent
Since we can consider Chern-Simons with BF case if $n$ is odd,
In three dimension, we can construct another model.
Let $E$ be a vector bundle with a Poisson structure on the fiber.
If we take $n=3$ in the equation (\ref{totbcsike}), 
we obtain $E\left[1\right]$.
A {\it P-structure} is defined 
on the graded vector bundle
$
E[1] \oplus T^*[2] M
$
by setting $n=3$ in the equation (\ref{csantiike}).
The abelian action is the $n=3$ case in the action (\ref{s0csike}):
\begin{eqnarray}
S_0 = && \int_{X} 
- \bb_{2 i} d \bphiike{}^{i}
+ \frac{k^{ab}}{2} \ba_{1 a} d \ba_{1 b},
\label{cs03ike}
\end{eqnarray}
and deformation is obtained as follows:
\begin{eqnarray}
S = && S_0 + g S_1, \nonumber \\ 
S_1 = && \int_{X} \left(
f_{1a}{}^i (\bphiike) \ba^a \bb_i 
+  \frac{1}{6} f_{2abc} (\bphiike) \ba^a \ba^b \ba^c
\right),
\label{cs1ike}
\end{eqnarray}
where we rewrite two structure functions $f_{1a}{}^i = g F_{11, a}{}^i$
and $\frac{1}{6} f_{2abc} = g F_{30, abc}$ for clarity.
If we substitute (\ref{cs1ike}) to the condition (\ref{s1s1ike}),
we obtain the identities on the structure functions 
$f_{1a}{}^i$ and $f_{2abc}$ as
\begin{eqnarray}
&& k^{ab} f_{1a}{}^i f_{1b}{}^j = 0, \nonumber \\ 
&& \left( \frac{\ldike}{\partial \bphiike^j} f_{1b}{}^i \right) 
f_{1c}{}^j
- \left( \frac{\ldike}{\partial \bphiike^j} f_{1c}{}^i \right) 
f_{1b}{}^j
+ k^{ef} f_{1e}{}^i f_{2fbc} = 0, \nonumber \\
&& \left( f_{1d}{}^j \frac{\ldike }{\partial \bphiike^j} f_{2abc}
- f_{1c}{}^j \frac{\ldike}{\partial \bphiike^j}
f_{2dab}
+ f_{1b}{}^j \frac{\ldike }{\partial \bphiike^j} f_{2cda}
- f_{1a}{}^j \frac{\ldike }{\partial \bphiike^j} f_{2bcd}
\right) \nonumber \\
&& 
+ k^{ef} (f_{2eab} f_{2cdf} 
+ f_{2eac} f_{2dbf} 
+ f_{2ead} f_{2bcf})
= 0.
\label{csideike}
\end{eqnarray}

The identities (\ref{csideike}) define a Courant algebroid.
In the definition of the Courant algebroid,
$\calMike$ is the space of a map $\bphiike$ from $\Pi T X$ to $M$
and $\calEike$ is the space of sections of 
$\Pi T^* X \oplus \bphiike^*(E[1])$.
We define a graded symmetric bilinear form 
$\bracketike{\cdot}{\cdot}$, a bilinear operation $\circ$, 
a bundle map $\rho$ and ${\cal D}$ by the antibracket with
the same as the equation (\ref{corresbase2ike}),
%
%
Then we can confirm that the classical master equation 
$\sbvike{S}{S} =0$, which derive 
the identity (\ref{csideike}) on structure functions $f$'s, 
is equivalent to the conditions $1$ to $5$ of a Courant algebroid 
(\ref{courantdefike}).
\cite{I3ike}
We take the basis of the section of the fiber $e^a = \ba_{1}{}^a$.
We calculate the operations 
$\circ$ and $\rho$ on the basis as follows:
\begin{eqnarray}
&& \ba^a \circ \ba^b = - k^{ac} k^{bd} f_{2cde} (\bphiike) \ba^e, 
\nonumber \\
&&\bracketike{\ba^a}{\ba^b} = k^{ab},
\nonumber \\
&& \rho(\ba^a) \bphiike^i = - f_{1c}{}^i (\bphiike) k^{ac}.
  \label{cscourant}
\end{eqnarray}
All deformations of a topological sigma model on the space $E[1]$ have
a Courant algebroid structure.


\subsection{General $n$}
\noindent
In $n$ dimensions,
the equation (\ref{s1s1ike}), $\sbvike{S_1}{S_1} = 0$, impose 
an algebroid structure on the space (\ref{totbundleike}), 
$\sum_{p=1}^{\lfloor \frac{n-1}{2} \rfloor}
E_p[p] \oplus E_p^*[n-p-1]$.
The algebroid structure is derived from 
the Batalin-Vilkovisky structure $\sbvike{S_1}{S_1} = 0$
of nonlinear gauge theories.
Now we obtain an infinite series of algebroids labeled by $n$.
We call this algebroid an {\it $n$-algebroid}.

In the previous section, we have found that 
the $n=2$ case defines a Lie algebroid on $T^*M$ and
the $n=3$ case defines a Courant algebroid on $E \oplus E^*$.
For $n \geq 4$ cases, we can easily calculate algebraic relations but
characterization of algebroid structures is still unknown.
Higher order generalization has also been discussed in \cite{Sevike}.

In the Chern-Simons case,
the equation (\ref{s1s1ike}), $\sbvike{S_1}{S_1} = 0$ impose 
an algebroid structure on the space (\ref{totbcsike}), 
$
\left(\sum_{p=1}^{\frac{n-3}{2}}
E_p[p] \oplus E_p^*[n-p-1] \right) \oplus E\left[\frac{n-1}{2}\right]
$.
In the previous section, we have found that 
the $n=3$ case defines a Courant algebroid on $E$.

\section{Quantum Version of Deformation}
\noindent
In the previous sections, we have considered a classical BV 
structure.
In this section, we discuss quantum version of deformation of 
a Batalin-Vilkovisky structure.
In this section, we discuss the BF case. 
We can make a similar discussion in the Chern-Simons with BF case.

In order to quantize a gauge theory, we must fix the gauge.
Gauge fixing is carried out by adding a gauge fixing term $S_{GF}$
to the classical action $S$.\cite{GPSike}
The gauge fixed quantum action is $S_q = S + S_{GF}$.

We need the {\it BV Laplacian}.
The BV Laplacian is defined as follows:
\begin{eqnarray}
\bDeltaike  F \equiv 
\sum_{p=0}^{\lfloor \frac{n-1}{2} \rfloor}
\frac{\ldike}{\partial \ba_p{}^{a_p}} 
\frac{\ldike }{\partial \bb_{n-p-1,a_p}} F
\label{bvlaplacianike}
\end{eqnarray}
%
The BV Laplacian satisfies the following identity:
\begin{eqnarray}
\bDeltaike(F \cdot G) = 
(\bDeltaike F) G + (-1)^{(n+1)|F|} \sbvike{F}{G} 
+ (-1)^{|F|} F \bDeltaike G,
\end{eqnarray}

In order for the generating functional to be gauge invariant
in the quantum sense, 
the following quantum master equation is required: 
\begin{eqnarray}
\sbvike{S_q}{S_q} - 2 i \hbar \bDeltaike S_q = 0,
\label{qmastike}
\end{eqnarray}
for the quantum action $S_q$.
In our $n$-algebroid topological sigma model, 
two terms are independently satisfied, i.e.
$\bDeltaike S_q = 0$ and $\sbvike{S_q}{S_q} = 0$.

${\cal O}$ is called an observable 
if an operator ${\cal O}$ satisfies the following equation:
\begin{eqnarray}
\sbvike{S_q}{{\cal O}} - i \hbar \bDeltaike {\cal O}  = 0.
\label{obsbrsike}
\end{eqnarray}

Generally, there are two kinds of observables.
One is the integration of a local function $F$ 
on the boundary $\partial X$.
Let ${X_r} \subset \partial X$ be a $r$-cycle 
on the boundary $\partial X$.
If $F$ has the form degree $r$, the integration of 
$F$ on the $r$ cycle ${X_r}$:
\begin{eqnarray}
{\cal O} = \int_{X_r} F(\bPhiike)
\label{Wilobsike}
\end{eqnarray}
is nontrivial and satisfies
(\ref{obsbrsike}).

Another observable is constructed from $\ba_0{}^a$.
We consider a function $F$ of $\ba_0{}^a$ 
and restrict $F$ on the boundary,
$
{\cal O}_{F} \equiv F
(\ba_0{}^a)|_{\partial X}.
$
We can confirm that 
the form degree zero part ${\cal O}_{F}^{(0)}$ of ${\cal O}_{F}$ is an local
observable with ghost number zero on the boundary.

The generating functional is defined by the path integral as 
\begin{eqnarray}
Z[{\cal O}_k] =
\int \prod_{p=0}^{[\frac{n-1}{2}]} 
{\cal D}{\ba_p{}} {\cal D}{\bb_{n-p-1}} \ 
e^{\frac{i}{\hbar}( S_q + \sum_r J_r {\cal O}_r)},
\end{eqnarray}
where 
${\cal D}{\ba_p{}} {\cal D}{\bb_{n-p-1}}$ is a path integral measure and
$J_k$ are source fields and ${\cal O}_k$ are observables and
$\hbar = g$.%

We consider $n=2$ case.
Let $X$ be a two-dimensional disc.
Note that classical deformation derives a Poisson structure on $T^*M$ 
in $n=2$ case.
On the other hand, quantum deformation derives
the deformation quantization on a Poisson manifold $M$. \cite{Kike}
The correlation function of two local observables 
${\cal O}_f^{(0)}$ and ${\cal O}_g^{(0)}$ derives 
the Kontsevich's star product formula \cite{CFike}
on a Poisson manifold:
\begin{eqnarray}
f * g(x) = \int_{\phi(\infty)=x} 
{\cal D}{\bphiike} {\cal D}{\bb_{1}} \ 
{\cal O}^{(0)}_f(\bphiike(1)) {\cal O}^{(0)}_g(\bphiike(0))
e^{\frac{i}{\hbar} S_q},
\end{eqnarray}
where $\bphiike=\ba_0{}$ and $0, 1, \infty$ are 
three distinct points at the boundary $\partial X$.

If we calculate the same correlation function 
for $S_0$, we obtain the usual product of functions $f$ and $g$:
\begin{eqnarray}
f(x) g(x) = \int_{\phi(\infty)=x} 
{\cal D}{\bphiike} {\cal D}{\bb_{1}} \ 
{\cal O}^{(0)}_f(\bphiike(1)) {\cal O}^{(0)}_g(\bphiike(0))
e^{\frac{i}{\hbar} S_{0q}},
\end{eqnarray}
where $S_{0q} = S_0 + S_{GF}$.
Therefore quantum deformation in $n=2$ is equivalent to
the star deformation on $C^{\infty}(M)$.

We can generalize this discussion to higher orders.
Deformation $S_0 \rightarrow S$ derives
a generalization of the star deformation to higher
dimensions as follows:
\begin{eqnarray}
m_k[{\cal O}_1, {\cal O}_2, \cdots, {\cal O}_k]
= \int \prod_{p=0}^{[\frac{n-1}{2}]} 
{\cal D}{\ba_p{}} {\cal D}{\bb_{n-p-1}} \ 
{\cal O}_1 {\cal O}_2 \cdots {\cal O}_k
e^{\frac{i}{\hbar}S_q},
\label{genstarike}
\end{eqnarray}
under the appropriate regularization and the boundary conditions,
where $S$ is the deformation (\ref{lapSike}) of the abelian 
topological sigma model
and ${\cal O}_r$'s are two kinds of observables at the boundary. 
%
%
%
The correlation functions satisfy the Ward-Takahashi identity 
derived from the gauge symmetry:
\begin{eqnarray}
\int \prod_{p=0}^{[\frac{n-1}{2}]} 
{\cal D}{\ba_p{}} {\cal D}{\bb_{n-p-1}} \ 
\bDeltaike \left({\cal O} e^{\frac{i}{\hbar}S_q} \right) = 0,
\end{eqnarray}
which leads a quantum geometric structure on the space of 
correlation functions.

\section{Summary and Outlook}
\noindent
We have discussed deformation of Batalin-Vilkovisky structures
of topological sigma models in $n$ dimensions.
We have constructed general theory of 
most general deformation in general $n$ dimensions.
We have analyzed structures in the case of $n=2$ and $3$ in detail.
In $n=2$ deformation of a BV structure produces a Lie algebroid structure,
and in $n=3$, deformation produces a Courant algebroid structure.
For $n \geq 4$, 
characterization of $n$-algebroids obtained by deformation of topological 
sigma models is still unknown and an open problem.

We have also discussed quantum version of deformation.
For $n=2$ case, the deformation on the disc $X$ 
is equivalent to the deformation quantization 
on a Poisson manifold $M$.

In $n=2$, there are two special important cases of deformations.
They are $A$-model and $B$-model. \cite{W2ike}
There are many investigations to analyze quantum moduli.
For reviews, \cite{BCOVike}, \cite{HKKPTVVZike} 
and references therein.

$n = 3$ quantum deformation is analyzed in \cite{HMike}.
For general $n \geq 4$, quantum structures are unknown.
If we analyze higher $n$ cases,
we will obtain interesting mathematical and physical structures.

In this article, 
we assume the $p$ and $n-p-1$ are nonnegative integers
in $E [p] \oplus E^*[n-p-1]$, where
we identify the $p=0$ bundle with 
a cotangent bundle $T^*[n-1]M$.
We will be able to 
generalize our discussions to negative integers $p$ and $n-p-1$.
A special case has been analyzed in \cite{BM2ike}\cite{II2ike}.

We need make analysis of all moduli, i.e. we should consider 
Kodaira-Spencer theory of Batalin-Vilkovisky structures.

%
%
%
%
%
%

%
%



\printindex
\end{document}